\documentclass[12pt]{article}
\usepackage{amsmath}
\usepackage{amssymb}
\begin{document}

\title{Gr\"uss and Gr\"uss-Voronovskaya-type estimates for some Bernstein-type polynomials of real and complex variables}
\author{Sorin G. Gal \\
University of Oradea\\
Department of Mathematics and Computer Science\\
Str. Universitatii Nr. 1\\
410087 Oradea, Romania\\
e-mail: galso@uoradea.ro \\
and \\
Heiner Gonska\\
University of Duisburg-Essen\\
Department of Mathematics\\
D-47048 Duisburg, Germany\\
e-mail: heiner.gonska@uni-due.de}
\date{}
\maketitle

{\bf Abstract.} The first aim of this paper is to prove a Gr\"uss-Voronovskaya estimate for Bernstein and for a class of Bernstein-Durrmeyer polynomials on $[0, 1]$.
Then, Gr\"uss and Gr\"uss-Voronovskaya estimates for their corresponding operators of complex variable on compact disks are obtained. Finally, the results are extended to Bernstein-Faber polynomials attached to compact sets in the complex plane.



\section{Introduction}

For $f:[0, 1]\to \mathbb{C}$, the complex Bernstein polynomials are defined by
$B_{n}(f)(z)=\sum_{k=0}^{n}{n \choose k}z^{k}(1-z)^{n-k}f(k/n)$, $z\in \mathbb{C}$, $n\in \mathbb{N}$.

A whole chapter in the book by DeVore and Lorentz \cite{DevLor} is devoted to the case of real-valued functions $f$, as is the case in thousands of papers dealing with the operator $B_n$ and many of its modifications. The recent book \cite{Gal} extensively treats the case of complex-valued functions and the phenomenon of overconvergence occuring in this case.

An early and very important result is the 1932 theorem of Voronovskaya stating that for $f \in C^2 [0,1]$ one has the uniform convergence
$$
\lim_{n \to \infty} n \cdot [B_n (f;x)-f(x)] = \frac{x(1-x)}{2} \cdot f''(x) , x \in [0,1].
$$
A quantitative from of it can be found in a hardly known booklet by Videnskij \cite{Vid} who showed that
$$
\left |n \cdot [B_n (f;x) - f(x)] - \frac{x(1-x)}{2} f''(x)\right | \le x(1-x) \cdot \omega_1 \left( f'',\sqrt{\frac{2}{n}}\right)
$$
where $\omega_1$ is the first order modulus of continuity. Refined forms of the latter inequality will be used in this note.

The second classical result guiding us in our present research is the well-known Gr\"uss inequality for positive linear functionals $L : C[0,1] \to \mathbb{R}$. This inequality gives an upper bound for the generalized Chebyshev functional
$$
T(f,g) := L(f\cdot g) - L(f) \cdot L(g) , \quad f,g \in C[0,1].
$$
For positive linear operators $H : C[0,1] \to C[0,1]$ reproducing constant functions and $x \in [0,1]$ the functional $L = \epsilon_x \circ H$, hence $L(f) = H(f;x)$, was investigated in this context for the first time in \cite{AGR}; several papers by other authors followed this approach. In the recent note \cite {GoRaRu} it was shown that
$$
|H(fg;x) - H(f;x) \cdot H(g;x)|
$$
$$
\le \frac14 \cdot \tilde{\omega}_1 (f; 2 \cdot \sqrt{H(e_2;x)-H(e_1;x)^2}) \cdot
\tilde{\omega}_1 (g;2\cdot \sqrt{H(e_2;x) - H(e_1;x)^2})
$$
where $\tilde{\omega}_1$ is the least concave majorant of $\omega_1$ and $e_i (x) = x^i$ for $x \in [0,1]$. For Bernstein operators we thus obtained
$$
|B_n  (fg;x) - B_n (f;x) \cdot B_n (g;x)| \le \frac14 \cdot \tilde{\omega}_1 \left( f; 2 \cdot \sqrt{\frac{x(1-x)}{n}} \right) \cdot \tilde{\omega}_1 \left( g;2 \cdot \sqrt{\frac{x(1-x)}{n}} \right) .
$$
The aim of the present note is two-fold. We extend the Gr\"uss-type inequality to complex Bernstein and related polynomials by obtaining both upper and lower estimates for functions being analytic in open discs centered at $0$ and having radius $R > 1$. Results for Bernstein-Faber operators are also included.

Our second goal is to refine Gr\"uss-type inequalities for Bernstein, P\u alt\u anea and Bernstein-Faber operators in the spirit of Voronovskaya's theorem by identifying, for example,
$$
\lim_{n \to \infty} n \cdot [B_n (fg;x) - B_n (f;x) \cdot B_n (g;x)] = x(1-x) f'(x) \cdot g'(x)
$$
in both the real and the complex case. Moreover, estimates for the corresponding differences will be given. The approach to give such Gr\"uss-Voronovskaya-type estimates appears to be new.

The remainder of our present paper is organized as follows. Section 2 deals with the classical Bernstein operators in the real case, where an appropriate smoothing (interpolation) technique will be used. This technique is also employed in Section 3 dealing with a most interesting family of operators $U_n^\rho$, $\rho > 0$, introduced by P\u alt\u anea. Section 4 contains extensions to complex Bernstein polynomials of analytic functions; here both Gr\"uss- and Gr\"uss-Voronovskaya-type estimates are given. In Section 5 a parallel development for complex genuine Bernstein-Durrmeyer operators is presented. Our note is completed by Section 6 on Bernstein-Faber operators for which both types of inequalities are given.

\section{Gr\"uss-Voronovskaya estimate for Bernstein o\-pe\-rators on $[0, 1]$}

Supposing that $f, g\in C^{2}[0, 1]$, it is natural to ask for the limit
$$\lim_{n\to \infty}n[B_{n}(f g)(x)-B_{n}(f)(x)B_{n}(g)(x)].$$
Taking into account that by simple calculation we have
$$n[B_{n}(f g)(x)-B_{n}(f)(x)B_{n}(g)(x)]$$
$$=n\left \{B_{n}(f g)(x)-f(x)g(x)-\frac{x(1-x)}{2n}(f(x)g(x))^{\prime \prime}\right .$$
$$\left .-g(x)\left [B_{n}(f)(x)-f(x)-\frac{x(1-x)}{2 n}f^{\prime \prime}(x)\right ]\right .$$
$$\left . -B_{n}(f)(x)\left [B_{n}(g)(x)-g(x)-\frac{x(1-x)}{2 n}g^{\prime \prime}(x)\right ]\right .$$
$$\left . +\frac{x(1-x)}{n}f^{\prime}(x)g^{\prime}(x)+g^{\prime \prime}(x)\cdot \frac{x(1-x)}{2 n}[f(x)-B_{n}(f)(x)]\right \},$$
passing to the limit it easily follows
$$\lim_{n\to \infty}n[B_{n}(f g)(x)-B_{n}(f)(x)B_{n}(g)(x)]=x(1-x)f^{\prime}(x)g^{\prime}(x).$$
This suggests us to prove the following result called by us Gr\"uss-Voronovskaya-type estimate.

{\bf Theorem 2.1.} {\it If $f, g\in C^{2}[0, 1]$, then for all $x\in [0, 1]$ and $n\in \mathbb{N}$ we have
$$n\left |B_{n}(f g)(x)-B_{n}(f)(x)B_{n}(g)(x)-\frac{x(1-x)}{n}f^{\prime}(x)g^{\prime}(x)\right | $$
$$\le \frac{x(1-x)}{2}\left [\tilde{\omega}_{1}\left ((f g)^{\prime \prime} ; \frac{1}{3\sqrt{n}}\right )+\|g\|\cdot \tilde{\omega}_{1}\left (f^{\prime \prime} ; \frac{1}{3\sqrt{n}}\right )+\|f\|\cdot \cdot \tilde{\omega}_{1}\left (g^{\prime \prime} ; \frac{1}{3\sqrt{n}}\right )\right .$$
$$\left .+\frac{1}{2 n}\|f^{\prime \prime}\|\cdot \|g^{\prime \prime}\|\right ].$$
where $\tilde{\omega}$ is the least concave majorant of the modulus of continuity and $\|\cdot \|$ denotes the uniform norm in $C[0, 1]$.}

{\bf Proof.} For $x\in [0, 1]$ and $n\in \mathbb{N}$ we have the decomposition formula
$$B_{n}(f g)(x)-B_{n}(f)(x)B_{n}(g)(x)-\frac{x(1-x)f^{\prime}(x)g^{\prime}(x)}{n}$$
$$=\left [B_{n}(f g)(x)-(f g)(x)-\frac{x(1-x)(f g)^{\prime \prime}(x)}{2 n}\right ]$$
$$-f(x)\left [B_{n}(g)(x)-g(x)-\frac{x(1-x)g^{\prime \prime}(x)}{2 n}\right ]$$
$$-g(x)\left [B_{n}(f)(x)-f(x)-\frac{x(1-x)f^{\prime \prime}(x)}{2 n}\right ]$$
$$+[g(x)-B_{n}(g)(x)]\cdot [B_{n}(f)(x)-f(x)].$$
Therefore, by using the quantitative estimate in \cite{Go1}, the Gr\"uss-Voronovskaya functional
$$G(B_{n}, f, g; x):=B_{n}(f g)(x)-B_{n}(f)(x)B_{n}(g)(x)-\frac{x(1-x)}{n}f^{\prime}(x)g^{\prime}(x)$$
has the upper bound
$$|G(B_{n}, f, g; x)|\le \frac{x(1-x)}{2 n}\tilde{\omega}_{1}\left ((f g)^{\prime \prime} ; \frac{1}{3\sqrt{n}}\right )$$
$$+|g(x)|\cdot \frac{x(1-x)}{2 n}\cdot \tilde{\omega}_{1}\left (f^{\prime \prime} ; \frac{1}{3\sqrt{n}}\right )+
|f(x)|\cdot \frac{x(1-x)}{2 n}\cdot \tilde{\omega}_{1}\left (g^{\prime \prime} ; \frac{1}{3\sqrt{n}}\right )$$
$$+|g(x)-B_{n}(g)(x)|\cdot |f(x)-B_{n}(f)(x)|.$$
On the other hand, by \cite{Pal}, p. 96, Corollary 4.1.9 we have
$$|f(x)-B_{n}(f)(x)|\le \frac{x(1-x)}{2 n}\|f^{\prime \prime}\|,$$
which replaced above leads to the estimate
$$|n G(B_{n}, f, g;x)|\le \frac{x(1-x)}{2}\tilde{\omega}_{1}\left ((f g)^{\prime \prime} ; \frac{1}{3\sqrt{n}}\right )
+\|g\|\cdot \frac{x(1-x)}{2}\cdot \tilde{\omega}_{1}\left (f^{\prime \prime} ; \frac{1}{3\sqrt{n}}\right )$$
$$+\|f\|\cdot \frac{x(1-x)}{2}\cdot \tilde{\omega}_{1}\left (g^{\prime \prime} ; \frac{1}{3\sqrt{n}}\right )
+\frac{x^{2}(1-x)^{2}}{4 n}\cdot \|f^{\prime \prime}\|\cdot \|g^{\prime \prime}\|$$
$$\le\frac{x(1-x)}{2}\left [\tilde{\omega}_{1}\left ((f g)^{\prime \prime} ; \frac{1}{3\sqrt{n}}\right )+\|g\|\cdot \tilde{\omega}_{1}\left (f^{\prime \prime} ; \frac{1}{3\sqrt{n}}\right )+\|f\|\cdot \cdot \tilde{\omega}_{1}\left (g^{\prime \prime} ; \frac{1}{3\sqrt{n}}\right )\right .$$
$$\left .+\frac{1}{2 n}\|f^{\prime \prime}\|\cdot \|g^{\prime \prime}\|\right ],$$
proving the theorem. $\hfill \square$

An immediate consequence of Theorem 2.1 is the following corollary.

{\bf Corollary 2.2.} {\it If $f, g\in C^{3}[0, 1]$ then for all $x\in [0, 1]$ and $n\in \mathbb{N}$ we have
$$n\left |B_{n}(f g)(x)-B_{n}(f)(x)B_{n}(g)(x)-\frac{x(1-x)}{n}f^{\prime}(x)g^{\prime}(x)\right |$$
$$\le C x(1-x)\cdot \frac{1}{\sqrt{n}}\left [\|(f g)^{\prime \prime \prime}\|+\|f^{\prime \prime}\|\cdot \|g^{\prime \prime}\|+
\|g\|\cdot \|f^{\prime \prime \prime}\|+\|f\|\cdot \|g^{\prime \prime \prime}\|\right ]$$
$$=O\left (\frac{x(1-x)}{\sqrt{n}}\right ),$$
where $C>0$ is an absolute constant (independent of $n$, $f$ and $g$).}

Also, we have :

{\bf Corollary 2.3.} {\it If $f, g\in C^{3}[0, 1]$ then for all $x\in [0, 1]$ and $n\in \mathbb{N}$ we have
$$n\left |B_{n}(f g)(x)-B_{n}(f)(x)B_{n}(g)(x)-\frac{x(1-x)}{n}f^{\prime}(x)g^{\prime}(x)\right |$$
$$\le C x(1-x)\cdot \frac{1}{\sqrt{n}}\cdot \max\{\|f\|, \|f^{\prime \prime \prime}\|\}\cdot \max\{\|g\|, \|g^{\prime \prime \prime}\|\},$$
where $C>0$ is an absolute constant (independent of $n$, $f$ and $g$).}

{\bf Proof.} The proof is immediate from Corollary 2.2 and taking into account the following remark in \cite{GoDis}, pp. 58-59 :

if $k\ge 2$ and $f\in C^{k}[a, b]$, then for any $0\le j\le k$ one has
$$\|f^{(j)}\|\le c\cdot \max\{\|f\|, \|f^{(k)}\|\},$$
with $c=c(k, j, b-a\}$, but independent of $f$. $\square \hfill$

{\bf Remarks.} 1) Since in the left hand-side of the estimate in Theorem 2.1 appear only $f^{\prime}$ and $g^{\prime}$, it is
natural to ask for an order of approximation when $f$ and $g$ are only in $C^{1}[0, 1]$.

2) If one function is constant and the other one is in $C^{2}[0, 1]$ or if both functions are linear, then while the left hand-side in the estimate in Theorem 2.1 is equal to zero, the estimate on the right hand-side is not zero, which shows that there is room for improvements.

As an answer to these remarks, we have the following result.

{\bf Theorem 2.4.} {\it Let $f, g\in C^{1}[0, 1]$ and $n\ge 1$. Then there is a constant $C$ independent of $n, f, g$ and $x$, such that
$$n\left |B_{n}(f g)(x)-B_{n}(f)(x)B_{n}(g)(x)-\frac{x(1-x)}{n}f^{\prime}(x)g^{\prime}(x)\right |$$
$$\le C x(1-x)\cdot \left \{\omega_{3}\left (f^{\prime};\frac{1}{n^{1/6}}\right )\cdot \omega_{3}\left (g^{\prime};\frac{1}{n^{1/6}}\right )\right .$$
$$\left .+\|f^{\prime}\|\cdot \omega_{3}\left (g^{\prime};\frac{1}{n^{1/6}}\right )+\|g^{\prime}\|\cdot \omega_{3}\left (f^{\prime};\frac{1}{n^{1/6}}\right )\right .$$
$$\left .+\max\left \{\frac{1}{n^{1/2}}\cdot \|f\|, \omega_{3}\left (f^{\prime};\frac{1}{n^{1/6}}\right )\right \}\cdot
\max\left \{\frac{1}{n^{1/2}}\cdot \|g\|, \omega_{3}\left (g^{\prime};\frac{1}{n^{1/6}}\right )\right \}\right \}.$$}

{\bf Proof.} In the considerations below, $C$ will always denote a constant independent of $n, f, g$ and $x$ , which may change its values during
the course of the proof.

For brevity, everywhere in this proof we will denote $G(f, g)(x)=B_{n}(f g)(x)-B_{n}(f)(x)B_{n}(g)(x)-\frac{x(1-x)}{n}f^{\prime}(x)g^{\prime}(x)$.
Since $G(f, g)$ is bilinear, for fixed $f, g\in C^{1}[0, 1]$ and $u, v\in C^{4}[0, 1]$ arbitrary, we can write
$$|G(f, g)(x)|=|G(f-u+u, g-v+v)(x)|$$
$$\le |G(f-u, g-v)(x)|+|G(u, g-v)(x)|+|G(f-u, v)(x)|+|G(u, v)(x)|.$$
Taking into account that by Theorem 4, pp. 854-855 in \cite{AGR} there exist $\eta, \theta\in [0, 1]$ with
$$B_{n}(fg)(x)-B_{n}(f)(x)\cdot B_{n}(g)(x)=f^{\prime}(\eta)\cdot g^{\prime}(\theta)\cdot \frac{x(1-x)}{n},$$
we easily get that
$$|nG(f, g)(x)|=|[f^{\prime}(\eta)\cdot g^{\prime}(\theta)-f^{\prime}(x)\cdot g^{\prime}(x)]x(1-x)|\le 2\|f^{\prime}\|\cdot \|g^{\prime}\|x(1-x).$$
Using this in the above estimate for $G(f, g)(x)$, we obtain
$$|G(f, g)(x)|$$
$$\le \frac{Cx(1-x)}{n}\cdot [\|(f-u)^{\prime}\|\cdot \|(g-v)^{\prime}\|+\|(f-u)^{\prime}\|\cdot \|v^{\prime}\|+\|u^{\prime}\|\cdot \|(g-v)^{\prime}\|]+|G(u, v)(x)|.$$
To estimate the term $G(u, v)(x)$, we use a simple consequence of Theorem 4 in \cite{GoRas}, stating that for $f\in C^{4}[0, 1]$ one has
$$\left |n[B_{n}(f)(x)-f(x)]-\frac{x(1-x)}{2}\cdot f^{\prime \prime}(x)\right |\le \frac{x(1-x)}{n}(\|f^{\prime \prime \prime}\|+\|f^{(4)}\|).$$
This fact allows us to write, for $u, v\in C^{4}[0, 1]$ and using the same decomposition as in the proof of Theorem 2.1 and the remark in \cite{GoDis}, pp. 58-59 (see the proof of Corollary 2.3),
$$|G(u, v)(x)|=\left |B_{n}(u v)(x)-B_{n}(u)(x)-B_{n}(v)(x)-\frac{x(1-x)}{n}\cdot u^{\prime}(x)\cdot v^{\prime}(x)\right |$$
$$\le \frac{x(1-x)}{n^{2}}\cdot (\|(u\cdot v)^{\prime \prime \prime}\|+\|(u\cdot v)^{(4)}\|)+\|u\|\cdot \frac{x(1-x)}{n^{2}}(\|v^{\prime \prime \prime}\|
+\|v^{(4)}\|)$$
$$+\|v\|\cdot \frac{x(1-x)}{n^{2}}(\|u^{\prime \prime \prime}\|+\|u^{(4)}\|)+\frac{x^{2}(1-x)^{2}}{n^{2}}\cdot \|u^{\prime \prime}\|\cdot \|v^{\prime \prime}\|$$
$$\le C\cdot \frac{x(1-x)}{n^{2}}\cdot \max\{\|u\|, \|u^{(4)}\|\}\cdot \max\{\|v\|, \|v^{(4)}\|\}.$$
Thus, we arrive at
$$|G(f, g)(x)|\le \frac{Cx(1-x)}{n}
\cdot \left [\|(f-u)^{\prime}\|\cdot \|(g-v)^{\prime}\|+\|(f-u)^{\prime}\|\cdot \|v^{\prime}\|\right .$$
$$\left .+\|u^{\prime}\|\cdot \|(g-v)^{\prime}\|
+\frac{1}{n}\cdot \max\{\|u\|, \|u^{(4)}\|\}\cdot \max\{\|v\|, \|v^{(4)}\|\}\right ].$$
Applying now to the right-hand side of the above inequality Lemma 3.1, p. 160 in \cite{GoLac}, for the particular cases $r=1$, $s=2$ there and choosing $u=f_{h, 3}$, $v=g_{h, 3}$, for all $h\in (0, 1]$ and $n\in \mathbb{N}$, it follows
$$|G(f, g)(x)|$$
$$\le \frac{Cx(1-x)}{n}\left \{\omega_{3}(f^{\prime};h)\cdot \omega_{3}(g^{\prime};h)+
\frac{1}{h}\cdot \omega_{1}(f;h)\cdot \omega_{3}(g^{\prime};h) + \omega_{3}(f^{\prime};h)\cdot \frac{1}{h}\omega_{1}(g ;h)\right .$$
$$\left .+\frac{1}{n}\cdot \max\left \{\|f\|, \frac{1}{h^{3}}\cdot \omega_{3}(f^{\prime};h)\right \}\cdot \max\left \{\|g\|, \frac{1}{h^{3}}\cdot \omega_{3}(g^{\prime};h)\right \}\right \}$$
$$\le \frac{Cx(1-x)}{n}\left \{\omega_{3}(f^{\prime};h)\cdot \omega_{3}(g^{\prime};h)+\|f^{\prime}\|\cdot \omega_{3}(g^{\prime};h)+\|g^{\prime}\|\cdot \omega_{3}(f^{\prime} ; h)\right .$$
$$\left .+\frac{1}{n}\cdot \max\left \{\|f\|, \frac{1}{h^{3}}\cdot \omega_{3}(f^{\prime};h)\right \}\cdot \max\left \{\|g\|, \frac{1}{h^{3}}\cdot \omega_{3}(g^{\prime};h)\right \}\right \}.$$
Choosing above $h=1/n^{1/6}$, we get
$$|G(f, g)(x)| \le \frac{Cx(1-x)}{n}$$
$$\cdot \left \{\omega_{3}(f^{\prime};1/n^{1/6})\cdot \omega_{3}(g^{\prime};1/n^{1/6})+\|f^{\prime}\|\cdot \omega_{3}(g^{\prime};1/n^{1/6})
+\|g^{\prime}\|\cdot \omega_{3}(f^{\prime} ; 1/n^{1/6})\right .$$
$$+\left .\max\left \{\frac{\|f\|}{n^{1/2}}, \omega_{3}(f^{\prime};1/n^{1/6})\right \}\cdot \max\left \{\frac{\|g\|}{n^{1/2}},  \omega_{3}(g^{\prime};1/n^{1/6})\right \}\right \}.$$
This proves the theorem. $\square \hfill$

In order to prove that Theorem 2.4 is best possible as for as order is concerned, we consider the following

{\bf Example 2.5.} Let $f(x)=e_{1}(x)=x$, $g(x)=e_{2}(x)=x^{2}$. Since $B_{n}(e_{0})(x)=1$ and $B_{n}(e_{1})(x)=x$, denoting $e_{k}(x)=x^{k}$,  by the recurrence formula in \cite{Andrica} (see also \cite{Gal}, p. 7)
$$B_{n}(e_{k+1})(x)=\frac{x(1-x)}{n}B_{n}^{\prime}(e_{k})(x)+x B_{n}(e_{k})(x),$$
we immediately obtain $B_{n}(e_{2})(x)=x^{2}+\frac{x(1-x)}{n}$ and
$$B_{n}(e_{3})(x)=x^{3}+\frac{3x^{2}(1-x)}{n}+\frac{x(1-x)(1-2x)}{n^{2}}.$$
Thus, we immediately obtain
$$G(B_{n}, f; g ; x)=B_{n}(e_{3})(x)-B_{n}(e_{1})(x)- B_{n}(e_{2})(x)\frac{x(1-x)}{n}e_{1}^{\prime}(x)\cdot e_{2}^{\prime}(x)$$
$$=x^{3}+\frac{3x^{2}(1-x)}{n}+\frac{x(1-x)(1-2x)}{n^{2}}-x^{3}-\frac{x^{2}(1-x)}{n} - \frac{2x^{2}(1-x)}{n}$$
$$=\frac{x(1-x)(1-2x)}{n^{2}}.$$
For $f$ and $g$ given as above, the estimate in Theorem 2.4 yields
$$n|G(B_{n}, e_{1}, e_{2} ;x)|\le Cx(1-x)\{0\cdot 0+ 1\cdot 0 + 2\cdot 0 +\max\{1/n^{1/2}, 0\}\cdot \max\{1/n^{1/2}, 0\}\}$$
$$=C\cdot \frac{x(1-x)}{n},$$
a fact which shows that the order in Theorem 2.4 cannot be improved in general.

\section{Results for P\u alt\u anea operators $U_{n}^{\rho}$ on $[0, 1]$}

In this section we extend the results from Section 2 to a class of one-parameter operators $U_{n}^{\rho}$ of Bernstein-Durrmeyer
type that preserve linear functions and constitute a link between the so-called "genuine Bernstein-Durrmeyer operators" $U_{n}$ and the classical
Bernstein operators $B_{n}$. The investigation on the operators in question started in a 2007 note by P\u alt\u anea \cite{Pal_2007}.
We recall some basic facts.

{\bf Definition 3.1.} Let $\rho >0$ and $n\in \mathbb{N}$. For $f\in C[0, 1]$ and $x\in [0, 1]$, let us define the polynomial operators
$$U_{n}^{\rho}(f)(x)=\sum_{k=1}^{n-1}\left (\int_{0}^{1}f(t)\mu^{\rho}_{n, k}(t)d t\right )p_{n, k}(x)+f(0)(1-x)^{n}+f(1)x^{n},$$
where $p_{n, k}(x)={n \choose k}x^{k}(1-x)^{n-k}$, $\mu_{n, k}^{\rho}(t):=\frac{t^{k \rho -1}(1-t)^{(n-k)\rho -1}}{B(k \rho, (n-k)\rho)}, 1\le k \le n-1$
and $B(x, y)=\int_{0}^{1}t^{x-1}(1-t)^{y-1} dt,$ $x, y >0$ is Euler's Beta function.

{\bf Remark.} For $\rho=1$ we obtain the genuine Bernstein-Durrmeyer operators given by
$$U_{n}(f)(x)=(n-1)\sum_{k=1}^{n-1}\left (\int_{0}^{1}f(t)p_{n-2, k-1}(t)d t\right )p_{n, k}(x)+(1-x)^{n}f(0)+x^{n}f(1),$$
while for each $f\in C[0, 1]$ we have $\lim_{\rho \to \infty}U_{n}^{\rho}(f)=B_{n}(f)$ uniformly.

Several properties of the operators $U_{n}^{\rho}$ were proved in the papers \cite{GoPal1}, \cite{GoPal2}.

The following quantitative Gr\"uss-Voronovskaya-type inequality holds.

{\bf Theorem 3.2.} {\it Let $f, g\in C^{1}[0, 1]$ and $n\ge 1$. Then there is a constant $C_{\rho}>0$ independent of $n, f, g$
and $x$, such that
$$n\left |U_{n}^{\rho}(f\cdot g)(x)-U_{n}^{\rho}(f)(x)\cdot U_{n}^{\rho}(g)(x)-\frac{\rho+1}{\rho n}\cdot x(1-x)f^{\prime}(x)\cdot g^{\prime}(x)\right |$$
$$\le C_{\rho}\cdot \frac{\rho+1}{\rho}\cdot x(1-x)\cdot \left \{\omega_{3}(f^{\prime} ; \delta_{n, \rho}^{1/6})\cdot \omega_{3}(g^{\prime} ; \delta_{n, \rho}^{1/6})\right .$$
$$+\left .\|f^{\prime}\|\cdot \omega_{3}(g^{\prime} ; \delta_{n, \rho}^{1/6})+\|g^{\prime}\|\cdot \omega_{3}(f^{\prime} ; \delta_{n, \rho}^{1/6})\right .$$
$$\left .+\max\{\delta_{n, \rho}^{1/2}\cdot \|f\|, \omega_{3}(f^{\prime} ; \delta_{n, \rho}^{1/6})\}\cdot \max\{\delta_{n, \rho}^{1/2}\cdot \|g\|, \omega_{3}(g^{\prime} ; \delta_{n, \rho}^{1/6})\}\right \},$$
where $\delta_{n, \rho}=\frac{\rho+1}{n\rho +1}$. For $\rho\to \infty$, the constants $C_{\rho}$ remain bounded.}

{\bf Proof.} Let $f, g\in C^{1}[0, 1]$. Applying Theorem 4, pp. 854-855 in \cite{AGR} as in the proof of Theorem 2.4, we can write
$$|G(U_{n}^{\rho}, f, g ;x)|:=|U_{n}^{\rho}(f g)(x)-U_{n}^{\rho}(f)(x)\cdot U_{n}^{\rho}(g)(x)-\frac{\rho +1}{\rho n}\cdot x(1-x)f^{\prime}(x)g^{\prime}(x)|$$
$$=|f^{\prime}(\eta)\cdot g^{\prime}(\theta)\cdot \frac{\rho +1}{\rho n +1}\cdot x(1-x)-\frac{\rho +1}{\rho n}\cdot x(1-x)\cdot f^{\prime}(x)\cdot g^{\prime}(x)|.$$
$$\le 2\frac{\rho+1}{\rho n}\cdot x(1-x)\cdot \|f^{\prime}\|\cdot \|g^{\prime}\|.$$
Reasoning as in the proof of Theorem 2.4, for $u, v\in C^{4}[0, 1]$ arbitrary, we can write
$$|G(U_{n}^{\rho}, f, g;x)|\le |G(f-u, g-v ;x)|+|G(u, g-v;x)|+|G(f-u, v ;x)|+|G(u, v ;x)|$$
$$\le 2\cdot \frac{\rho+1}{n\rho}x(1-x)[\|(f-u)^{\prime}\|\cdot \|(g-v)^{\prime}\|+\|(f-u)^{\prime}\|\cdot \|v^{\prime}\|+\|u^{\prime}\|
\cdot \|(g-v)^{\prime}\|]$$
$$+|G(u, v)(x)|,$$
where
$$|G(u, v)(x)|=|U_{n}^{\rho}(u v)(x)-U_{n}^{\rho}(u)(x)\cdot U_{n}^{\rho}(v)(x)-\frac{\rho +1}{\rho n}\cdot x(1-x)u^{\prime}(x)\cdot v^{\prime}(x)|$$
$$=|U_{n}^{\rho}(u v)(x)-(u v)(x)-\frac{\rho +1}{2\rho n}\cdot x(1-x)(u v)^{\prime \prime}(x)$$
$$-u(x)[U_{n}^{\rho}(v)(x)-v(x)-\frac{\rho +1}{2\rho n}\cdot x(1-x)\cdot v^{\prime \prime}(x)]$$
$$-v(x)[U_{n}^{\rho}(u)(x)-u(x)-\frac{\rho +1}{2\rho n}\cdot x(1-x)\cdot u^{\prime \prime}(x)]+[v(x)-U_{n}^{\rho}(v)(x)]\cdot [U_{n}^{\rho}(u)(x)-u(x)].$$

Corollary 5.1 in \cite{GoPal1} gives for $f \in C^2 [0,1]$
$$
|U_n^\rho (f;x) - f(x)- \frac12 U_n^\rho ((e_1-x)^2;x) \cdot f''(x)|
$$
$$
\le U_n^\rho ((e_1-x)^2;x) \left\{ \frac56 \cdot \sqrt{A} \cdot \omega_1 (f'';\sqrt{B}) + \frac{13}{16} \cdot \omega_2 (f'';\sqrt{B})\right\}
$$
where
$$
A := \frac{[M_3(x)]^2}{M_2(x) M_4(x)}, \quad B := \frac{M_4(x)}{M_2(x)} , \quad M_r (x) = U_n^\rho ((e_1-x)^r;x).
$$
Explicitly,
$$
A = \frac{(\rho + 2)^2(X')^2 (n\rho + 3)}{(n\rho + 2) \{[3\rho (\rho+1)n - 6 (\rho^2 + 3 \rho + 3] X + (\rho + 2)(\rho + 3)\}} ,
$$
$$
B = \frac{3[\rho (\rho + 1)n - 2 (\rho^2 + 3\rho + 3)] X + (\rho + 2)(\rho + 3)}{(n\rho + 2)(n\rho + 3)} ,
$$
where $X := x (1-x)$.

Moreover, $U_n^\rho ((e_1 - x)^2;x) = \frac{(\rho + 1)x(1-x)}{n\rho + 1}$.

We slightly modify the Voronovskaya expression from above and consider
$$
\left| U_n^\rho (f;x) - f(x) - \frac{\rho+1}{2\rho n} x (1-x) \cdot f''(x) \right|
$$
$$
\le \left| U_n^\rho (f;x) - f(x) - \frac{\rho+1}{2(\rho n  + 1)} \cdot x (1-x) \cdot f''(x) \right|
$$
$$
 + \left| \frac{\rho + 1}{2(\rho n + 1)} - \frac{\rho + 1}{2\rho n} \right| x (1-x)| \cdot |f''(x)|
$$
$$
\le \frac{(\rho + 1)x(1-x)}{n\rho + 1} \left\{ \frac56 \cdot \sqrt{A} \cdot \omega_1 (f'';\sqrt{B}) + \frac{13}{16} \cdot \omega_2 (f''; \sqrt{B} \right\}
$$
$$
+ \frac{1+\rho}{2(\rho n + 1)\rho n} \cdot x(1-x) \cdot |f''(x)|
$$
$$
\le \frac{(\rho + 1)x(1-x)}{n\rho + 1} \cdot \left\{ \frac56 \sqrt{A\cdot B} \cdot ||f'''|| + \frac{13}{16} \cdot B \cdot ||f^{(4)}||\right\}
$$
$$
+ \frac{1+\rho}{2(\rho n + 1)\cdot \rho n} \cdot x (1-x) \cdot ||f''|| \;\;\; (f \in C^4 [0,1])
$$
$$
\le \frac{(\rho + 1) x (1-x)}{n\rho + 1} \cdot \left\{ \sqrt{A\cdot B} \cdot ||f'''|| + B \cdot ||f^{(4)}|| + \frac{1}{\rho \cdot n} \cdot ||f''||\right\}
$$
$$
\le C \cdot \frac{(\rho + 1) x (1-x)}{n\rho + 1} \cdot \frac{1}{n} \cdot \{ ||f''|| + ||f'''|| + ||f^{(4)}||\}
$$
with constants $C = C (\rho)$ that remain bounded for $n$ fixed and $\rho \to \infty$.

Collecting all these, they imply for $u,v \in C^4 [0,1]$
$$
|G(u,v)(x)| \le C \cdot \frac{(\rho + 1) x (1-x)}{n\rho + 1} \cdot \frac{1}{n} \cdot \{ ||(uv)''|| + ||(uv)'''|| + ||(uv)^{(4)}||\}
$$
$$
+ ||u|| \cdot C \cdot \frac{(\rho + 1) x (1-x)}{n\rho + 1} \cdot \frac{1}{n} \cdot \{ ||v''|| + ||v'''|| + ||v^{(4)}||\}
$$
$$
+ ||v|| \cdot C \cdot \frac{(\rho + 1) x (1-x)}{n\rho + 1} \cdot \frac{1}{n} \cdot \{||u''|| + ||u'''|| + ||u^{(4)}||\}
$$
$$
+ \left[ \frac{\rho + 1}{2(n\rho + 1)}\right]^2 \cdot x^2 (1-x)^2 \cdot ||u''|| \cdot ||v''||
$$
$$
\le C \cdot x (1-x) \cdot \frac{\rho + 1}{n\rho + 1} \cdot \frac{\rho + 1}{\rho} \cdot \frac{1}{n} \cdot \max \{||u||,||u^{(4)}||\} \cdot \max \{||v||,||v^{(4)}||\} .
$$
Collecting the above information we have now for
$
f,g \in C^1 [0,1]:
$
$$
|G(U_n^\rho, f,g,x)|
$$
$$
\le C \cdot \frac{x(1-x)}{n} \cdot \frac{\rho + 1}{\rho} \cdot [||(f-u)'|| \cdot ||(g-v)'||
$$
$$
+ ||(f-u)'|| \cdot ||v'|| + ||u'|| \cdot ||(g-v)'|| + \frac{\rho + 1}{n\rho + 1} \cdot \max \{||u|| , ||u^{(4)}||\}Ê\cdot \max \{||v||,||v^{(4)}||\}].
$$
The rest of the proof follows the pattern of that of Theorem 2.4 observing that there is the extra constant $\frac{\rho + 1}{\rho} \in (1,\infty)$ and that $h = \sqrt[6]{\frac{\rho + 1}{n \rho + 1}} \in (0,1)$ is an appropriate choice.

Multiplying both sides by $n$ then gives the desired result. \hfill $\square$

{\bf Remarks.} (i) For $\rho=1$ we get $\delta_{n, \rho}=\frac{2}{n+1}$ and a Gr\"uss-Voronovskaya inequality for the genuine Bernstein-Durrmeyer operators.

(ii) For $n$ fixed and $\rho\to \infty$, we recapture the result for Bernstein polynomials in Theorem 2.4, as the constants $C_{\rho}$ remain bounded for $\rho \to \infty$.

\section{Results for complex Bernstein polynomials}

In this section we extend the Gr\"uss and the Gr\"uss-Voronovskaya estimates for complex Bernstein polynomials attached to analytic functions in compact disks.

Firstly, the following Gr\"uss-type inequality holds.

{\bf Theorem 4.1.} {\it Suppose that $R>1$ and $f, g:\mathbb{D}_{R}\to \mathbb{C}$ are analytic in $\mathbb{D}_{R}=\{z\in \mathbb{C} ; |z|<R\}$,
that is $f(z)=\sum_{k=0}^{\infty}a_{k}z^{k}$ and $f(z)=\sum_{k=0}^{\infty}b_{k}z^{k}$ for all $z\in \mathbb{D}_{R}$.

Let $1\le r <R$. Denoting $\|f\|_{r}=\max\{|f(z)| ; |z|\le r\}$, for all $n\in \mathbb{N}$ we have
$$\|B_{n}(f g)-B_{n}(f)B_{n}(g)\|_{r}\le \frac{6(1+r)}{n}\sum_{m=0}^{\infty}m^{2}\left [\sum_{j=0}^{m}|a_{j}|\cdot |b_{m-j}|\right ]r^{m-1},$$
where $\sum_{m=0}^{\infty}m^{2}\left [\sum_{j=0}^{m}|a_{j}|\cdot |b_{m-j}|\right ]r^{m-1}<+\infty$.}

{\bf Proof.} Denote $e_{m}(z)=z^{m}$. Since $f(z)g(z)=\sum_{m=0}^{\infty}c_{m}z^{m}$, where $c_{m}=\sum_{j=0}^{m}a_{j}b_{m-j}$, it follows
$$B_{n}(f g)(z)=\sum_{m=0}^{\infty}\left [\sum_{j=0}^{m}a_{j}b_{m-j}\right ]B_{n}(e_{m})(z).$$
Also,
$$B_{n}(f)(z)=\sum_{k=0}^{m}a_{k}B_{n}(e_{k})(z),\,\,  B_{n}(g)(z)=\sum_{k=0}^{m}b_{k}B_{n}(e_{k})(z)$$
and
$$B_{n}(f)(z)B_{n}(g)(z)=\sum_{m=0}^{\infty}\left [\sum_{j=0}^{m}a_{j}b_{m-j}B_{n}(e_{j})(z)B_{n}(e_{m-j})(z)\right ],$$
which immediately implies
$$|B_{n}(f g)(z)-B_{n}(f)(z)B_{n}(g)(z)|$$
$$=\left |\sum_{m=0}^{\infty}\left [\sum_{j=0}^{m}a_{j}b_{m-j}\left (B_{n}(e_{m})(z)-B_{n}(e_{j})(z)B_{n}(e_{m-j})(z)\right )\right ]\right |$$
$$\le \sum_{m=0}^{\infty}\left [\sum_{j=0}^{m}|a_{j}|\cdot |b_{m-j}|\cdot |B_{n}(e_{m})(z)-B_{n}(e_{j})(z)B_{n}(e_{m-j})(z)|\right ].$$
Then, we get
$$|B_{n}(e_{m})(z)-B_{n}(e_{j})(z)B_{n}(e_{m-j})(z)|$$
$$\le |B_{n}(e_{m})(z)-e_{m}(z)|+|e_{j}(z)\cdot e_{m-j}(z)-B_{n}(e_{j})(z)B_{n}(e_{m-j})(z)|$$
$$\le |B_{n}(e_{m})(z)-e_{m}(z)|+|e_{j}(z)|\cdot |e_{m-j}(z)-B_{n}(e_{m-j})(z)|$$
$$+|B_{n}(e_{m-j})(z)|\cdot |e_{j}(z)-B_{n}(e_{j})(z)|.$$
Taking into account that for all $|z|\le r$, $n, k\in \mathbb {N}$, we have $|B_{n}(e_{k})(z)|\le r^{k}$ (see e.g. \cite{Lor}, relationship (4), pp. 88) and $|B_{n}(e_{k})(z)-e_{k}(z)|\le \frac{3r(1+r)}{2 n}k(k-1)r^{k-2}$ (see e.g. \cite{Gal}, p. 8), from the above inequality it easily follows
$$|B_{n}(e_{m})(z)-B_{n}(e_{j})(z)B_{n}(e_{m-j})(z)|$$
$$\le \frac{3r(1+r)}{2 n}r^{m-2}[m(m-1)+(m-j)(m-j-1)+j(j-1)]$$
$$=\frac{3r(1+r)}{n}r^{m-2}[m^{2}-m-m j+j^{2}]\le \frac{6(1+r)}{n}m^{2} r^{m-1},$$
which leads to the inequality
$$|B_{n}(f g)(z)-B_{n}(f)(z)B_{n}(g)(z)|\le \frac{6(1+r)}{n}\cdot \sum_{m=0}^{\infty}m^{2}\left [\sum_{j=0}^{m}|a_{j}|\cdot |b_{m-j}|\right ]r^{m-1}.$$
It remains to show that $\sum_{m=0}^{\infty}m^{2}\left [\sum_{j=0}^{m}|a_{j}|\cdot |b_{m-j}|\right ]r^{m-1}<\infty$.
Indeed, since $f$ and $g$ are analytic it follows that the series $f(z)=\sum_{k=0}^{\infty}a_{k}z^{k}$ and $f(z)=\sum_{k=0}^{\infty}b_{k}z^{k}$
converges uniformly for $1\le r<R$, that is the series $\sum_{k=0}^{\infty}|a_{k}|r^{k}$ and $\sum_{k=0}^{\infty}|b_{k}|r^{k}$ converge for all
$1\le r<R$. Then, by Mertens' theorem (see e.g. \cite{Rud}, p. 74, Theorem 3.50) their (Cauchy) product is a convergent series and therefore
$$\sum_{m=0}^{\infty}\left [\sum_{j=0}^{m}|a_{j}|\cdot |b_{m-j}|\right ]r^{m}$$
is a convergent series for all $1\le r<R$. Denoting $A_{m}=\sum_{j=0}^{m}|a_{j}|\cdot |b_{m-j}|$, this means that the power series
$F(z)=\sum_{m=0}^{\infty}A_{m}z^{m}$ is uniformly convergent for $|z|\le r$, which implies that $F^{\prime \prime}(z)=\sum_{m=0}^{\infty}m(m-1)A_{m}z^{m-2}$ also is uniformly convergent for $|z|\le r$, that is
$\sum_{m=0}^{\infty}m(m-1)A_{m}r^{m-2}<\infty$. This immediately implies that $\sum_{m=0}^{\infty}m^{2}\left [\sum_{j=0}^{m}|a_{j}|\cdot |b_{m-j}|\right ]r^{m-1}<\infty$. $\hfill \square$

The Gr\"uss-Voronovskaya-type estimate follows.

{\bf Theorem 4.2.} {\it Suppose that $R > r \ge 1$ and $f, g:\mathbb{D}_{R}\to \mathbb{C}$ are analytic in $\mathbb{D}_{R}=\{z\in \mathbb{C} ; |z|<R\}$, that is $f(z)=\sum_{k=0}^{\infty}a_{k}z^{k}$ and $f(z)=\sum_{k=0}^{\infty}b_{k}z^{k}$ for all $z\in \mathbb{D}_{R}$.

Then, for all $n\in \mathbb{N}$ and $|z|\le r$ we have
$$\left |B_{n}(f g)(z)-B_{n}(f)(z)B_{n}(g)(z)-\frac{z(1-z)f^{\prime}(z)g^{\prime}(z)}{n}\right |\le \frac{C(r, f, g)}{n^{2}},$$
with $C(r, f, g)$ independent of $n$ and depending on $r, f, g$.}

{\bf Proof.} By the decomposition in the proof of Theorem 2.1, we get
$$B_{n}(f g)(z)-B_{n}(f)(z)B_{n}(g)(z)-\frac{z(1-z)f^{\prime}(z)g^{\prime}(z)}{n}$$
$$=\left [B_{n}(f g)(z)-(f g)(z)-\frac{z(1-z)(f g)^{\prime \prime}(z)}{2 n}\right ]$$
$$-f(z)\left [B_{n}(g)(z)-g(z)-\frac{z(1-z)g^{\prime \prime}(z)}{2 n}\right ]$$
$$-g(z)\left [B_{n}(f)(z)-f(z)-\frac{z(1-z)f^{\prime \prime}(z)}{2 n}\right ]$$
$$+[g(z)-B_{n}(g)(z)]\cdot [B_{n}(f)(z)-f(z)].$$
Passing to modulus with $|z|\le r$ and taking into account the estimates in Theorems 1.1.2 and 1.1.3 in \cite{Gal}, p. 6 and p. 9, we get
$$\left |B_{n}(f g)(z)-B_{n}(f)(z)B_{n}(g)(z)-\frac{z(1-z)f^{\prime}(z)g^{\prime}(z)}{n}\right |$$
$$\le \left |B_{n}(f g)(z)-(f g)(z)-\frac{z(1-z)(f g)^{\prime \prime}(z)}{2 n}\right |$$
$$+|f(z)|\left |B_{n}(g)(z)-g(z)-\frac{z(1-z)g^{\prime \prime}(z)}{2 n}\right |$$
$$+|g(z)|\left |B_{n}(f)(z)-f(z)-\frac{z(1-z)f^{\prime \prime}(z)}{2 n}\right |
+|g(z)-B_{n}(g)(z)|\cdot |B_{n}(f)(z)-f(z)|$$
$$\le \frac{C_{1}(r, f, g)}{n^{2}}+\|f\|_{r}\cdot \frac{C_{2}(r, g)}{n^{2}}+\|g\|_{r}\cdot \frac{C_{3}(r, f)}{n^{2}}
+\frac{C_{4}(r, g)}{n}\cdot \frac{C_{5}(r, f)}{n}\le \frac{C(r, f, g)}{n^{2}},$$
for all $n\in \mathbb{N}$ and $|z|\le r$, with $C(r, f, g)>0$ independent of $n$ and depending on $r, f, g$. $\hfill \square$

In what follows, the above theorem is used to obtain a lower estimate.

{\bf Corollary 4.3.} {\it Suppose that $R > r \ge 1$ and $f, g:\mathbb{D}_{R}\to \mathbb{C}$ are analytic in $\mathbb{D}_{R}=\{z\in \mathbb{C} ; |z|<R\}$, that is $f(z)=\sum_{k=0}^{\infty}a_{k}z^{k}$ and $f(z)=\sum_{k=0}^{\infty}b_{k}z^{k}$ for all $z\in \mathbb{D}_{R}$. If $f$ and $g$ are not constant functions, then for any $1\le r <R$, there exists a constant $K(r, f, g)$ (depending on $r$, $f$ and $g$), such that
$$\|B_{n}(f g)-B_{n}(f)B_{n}(g)\|_{r}\ge \frac{K(r, f, g)}{n}, \, n\in \mathbb{N}.$$}

{\bf Proof.}  We can write
$$B_{n}(f g)(z)-B_{n}(f)(z)B_{n}(g)(z)=\frac{1}{n}\left \{z(1-z)f^{\prime}(z)g^{\prime}(z)+\frac{1}{n}\right .$$
$$\left .\left [n^{2}\left (B_{n}(f g)(z)-B_{n}(f)(z)B_{n}(g)(z)-\frac{z(1-z)}{n}f^{\prime}(z)g^{\prime}(z)\right )\right ]\right \}.$$
Applying to the above identity the obvious inequality
$$\|F+G\|_{r}\ge | \|F\|_{r}-\|G\|_{r} |\ge \|F\|_{r} - \|G\|_{r},$$
and denoting $e_{1}(z)=z$, we obtain
$$\|B_{n}(f g)-B_{n}(f)B_{n}(g)\|_{r}\ge \frac{1}{n}\left \{\|e_{1}(1-e_{1})f^{\prime}g^{\prime}\|_{r}-\frac{1}{n}\right .$$
$$\left .\left [n^{2}\left \|B_{n}(f g)-B_{n}(f)B_{n}(g)-\frac{e_{1}(1-e_{1})}{n}f^{\prime}g^{\prime}\right \|_{r}\right ]\right \}.$$
Since $f$ and $g$ are not constant functions, we get $\|e_{1}(1-e_{1})f^{\prime}g^{\prime}\|_{r}>0$. Indeed, supposing the contrary, it follows
that $z(1-z)f^{\prime}(z)g^{\prime}(z)=0$, for all $|z|\le r$, which by analyticity easily implies that $f^{\prime}(z)g^{\prime}(z)=0$, for all
$|z|\le r$. Since the zeroes of analytic functions are isolated, it easily follows that $f$ is a constant function or $g$ is a constant
function, on $|z|\le r$, contradicting the hypothesis.

Taking into account that by Theorem 4.2 we get
$$n^{2}\left \|B_{n}(f g)-B_{n}(f)B_{n}(g)-\frac{e_{1}(1-e_{1})}{n}f^{\prime}g^{\prime}\right \|_{r}\le K(r, f, g)$$
and that $\frac{1}{n}\to 0$, there exists an index $n_{0}$ (depending only on $r, f, g$), such that  for all $n\ge n_{0}$ we have
$$\|e_{1}(1-e_{1})f^{\prime}g^{\prime}\|_{r}-\frac{1}{n}\left [n^{2}\left \|B_{n}(f g)-B_{n}(f)B_{n}(g)-\frac{e_{1}(1-e_{1})}{n}f^{\prime}g^{\prime}\right \|_{r}\right ]$$
$$\ge \frac{\|e_{1}(1-e_{1})f^{\prime}g^{\prime}\|_{r}}{2}>0,$$
which for all $n\ge n_{0}$ implies
$$\|B_{n}(f g)-B_{n}(f)B_{n}(g)\|_{r}\ge \frac{1}{n}\cdot \frac{\|e_{1}(1-e_{1})f^{\prime}g^{\prime}\|_{r}}{2}.$$
For $1\le n <n_{0}$, we obviously have
$$\|B_{n}(f g)-B_{n}(f)B_{n}(g)\|_{r}\ge \frac{M(r, n, f, g)}{n},$$
with $M(r, n, f, g)=n\cdot \|B_{n}(f g)-B_{n}(f)B_{n}(g)\|_{r}$. Since if $f$ and $g$ are not constant function we have $\|B_{n}(f g)-B_{n}(f)B_{n}(g)\|_{r}>0$, for all $n\in \mathbb{N}$, finally we get
$$\|B_{n}(f g)-B_{n}(f)B_{n}(g)\|_{r}\ge \frac{K(r, f, g)}{n}, \, n\in \mathbb{N},$$
where $K(r, f, g)=\min\left \{M(r, 1, f, g), ..., M(r, n_{0}-1, f, g), \frac{\|e_{1}(1-e_{1})f^{\prime}g^{\prime}\|_{r}}{2}\right \}$.
$\hfill \square$
As an immediate consequence of Theorem 4.1 and Corollary 4.3, we obtain the following exact estimate.

{\bf Corollary 4.4.} {\it Suppose that $R > r \ge 1$ and $f, g:\mathbb{D}_{R}\to \mathbb{C}$ are analytic in $\mathbb{D}_{R}=\{z\in \mathbb{C} ; |z|<R\}$, that is $f(z)=\sum_{k=0}^{\infty}a_{k}z^{k}$ and $f(z)=\sum_{k=0}^{\infty}b_{k}z^{k}$ for all $z\in \mathbb{D}_{R}$. If $f$ and $g$ are not constant functions, then for any $1\le r <R$ we have
$$\|B_{n}(f g)-B_{n}(f)B_{n}(g)\|_{r}\sim \frac{1}{n}, \, n\in \mathbb{N},$$
where the constants in the equivalence are independent of $n$ but depend on $r, f, g$.}

\section{Results for complex genuine Bernstein-Durr\-me\-yer operators}

The results in the previous section can be extended for the complex operators $U_{n}^{\rho}(f)(z)$
with arbitrary $\rho > 0$, but for simplicity of calculation, we consider here only results in the particular case
$\rho=1$, when $U_{n}^{\rho}$ reduce to the genuine Bernstein-Durrmeyer polynomials $U_{n}$, defined as in  Remark after Definition 3.1.

Firstly, the following Gr\"uss-type inequality holds.

{\bf Theorem 5.1.} {\it Suppose that $R>1$ and $f, g:\mathbb{D}_{R}\to \mathbb{C}$ are analytic in $\mathbb{D}_{R}=\{z\in \mathbb{C} ; |z|<R\}$,
that is $f(z)=\sum_{k=0}^{\infty}a_{k}z^{k}$ and $f(z)=\sum_{k=0}^{\infty}b_{k}z^{k}$ for all $z\in \mathbb{D}_{R}$.

Let $1\le r <R$. Denoting $\|f\|_{r}=\max\{|f(z)| ; |z|\le r\}$, for all $n\in \mathbb{N}$ we have
$$\|U_{n}(f g)-U_{n}(f)U_{n}(g)\|_{r}\le \frac{4}{n}\sum_{m=0}^{\infty}m^{2}\left [\sum_{j=0}^{m}|a_{j}|\cdot |b_{m-j}|\right ]r^{m},$$
where $\sum_{m=0}^{\infty}m^{2}\left [\sum_{j=0}^{m}|a_{j}|\cdot |b_{m-j}|\right ]r^{m}<+\infty$.}

{\bf Proof.} Denote $e_{m}(z)=z^{m}$. Similar to the proof of Theorem 4.1 we get
$$|U_{n}(e_{m})(z)-U_{n}(e_{j})(z)U_{n}(e_{m-j})(z)|$$
$$\le |U_{n}(e_{m})(z)-e_{m}(z)|+|e_{j}(z)|\cdot |e_{m-j}(z)-U_{n}(e_{m-j})(z)|$$
$$+|U_{n}(e_{m-j})(z)|\cdot |e_{j}(z)-U_{n}(e_{j})(z)|.$$
Taking into account that for all $|z|\le r$, $n, k\in \mathbb {N}$, we have $|U_{n}(e_{k})(z)|\le r^{k}$ (see \cite{Gal3}, Corollary 2.3, (i)) and $|U_{n}(e_{k})(z)-e_{k}(z)|\le \frac{2k(k-1)}{n}r^{k}$ (see \cite{Gal3}, p. 1916), from the above inequality it easily follows
$$|U_{n}(e_{m})(z)-U_{n}(e_{j})(z)U_{n}(e_{m-j})(z)|\le \frac{2 r^{m}}{n}[m(m-1)+(m-j)(m-j-1)+j(j-1)]$$
$$\le \frac{4 m^{2} r^{m}}{n}$$
and following the lines in the proof of Theorem 4.1 we obtain
$$|U_{n}(f g)(z)-U_{n}(f)(z)U_{n}(g)(z)|\le \frac{4}{n}\cdot \sum_{m=0}^{\infty}m^{2}\left [\sum_{j=0}^{m}|a_{j}|\cdot |b_{m-j}|\right ]r^{m}.$$
The theorem is proved. $\hfill \square$.

The Gr\"uss-Voronovskaya-type estimate follows.

{\bf Theorem 5.2.} {\it Suppose that $R > r \ge 1$ and $f, g:\mathbb{D}_{R}\to \mathbb{C}$ are analytic in $\mathbb{D}_{R}=\{z\in \mathbb{C} ; |z|<R\}$, that is $f(z)=\sum_{k=0}^{\infty}a_{k}z^{k}$ and $f(z)=\sum_{k=0}^{\infty}b_{k}z^{k}$ for all $z\in \mathbb{D}_{R}$.

Then, for all $n\in \mathbb{N}$ and $|z|\le r$ we have
$$\left |U_{n}(f g)(z)-U_{n}(f)(z)U_{n}(g)(z)-\frac{2z(1-z)f^{\prime}(z)g^{\prime}(z)}{n}\right |\le \frac{C(r, f, g)}{n^{2}},$$
with $C(r, f, g)$ independent of $n$ and depending on $r, f, g$.}

{\bf Proof.} By Theorem 2.4 in \cite{Gal3}, for all $|z|\le r$ and $n\in \mathbb{N}$, we have
$$\left |U_{n}(f)(z)-f(z)-\frac{z(1-z)f^{\prime \prime}(z)}{n+1}\right |\le \frac{M_{r}(f)}{n^{2}},$$
which immediately implies
$$\left |U_{n}(f)(z)-f(z)-\frac{z(1-z)}{n}\cdot f^{\prime \prime}(z)\right |$$
$$\le |U_{n}(f)(z)-f(z)-\frac{z(1-z)}{n+1}\cdot f^{\prime \prime}(z)|+|z||1-z|\cdot |f^{\prime \prime}(z)|\cdot \left |\frac{1}{n+1}-\frac{1}{n}\right |\le \frac{C_{r}(f)}{n^{2}}.$$
But we have the decomposition
$$U_{n}(f g)(z)-U_{n}(f)(z)U_{n}(g)(z)-\frac{2z(1-z)f^{\prime}(z)g^{\prime}(z)}{n}$$
$$=\left [U_{n}(f g)(z)-(f g)(z)-\frac{z(1-z)(f g)^{\prime \prime}(z)}{n}\right ]$$
$$-f(z)\left [U_{n}(g)(z)-g(z)-\frac{z(1-z)g^{\prime \prime}(z)}{n}\right ]$$
$$-g(z)\left [U_{n}(f)(z)-f(z)-\frac{z(1-z)f^{\prime \prime}(z)}{n}\right ]$$
$$+[g(z)-U_{n}(g)(z)]\cdot [U_{n}(f)(z)-f(z)].$$
Passing to modulus with $|z|\le r$ and taking into account the above estimate and that in Theorem 5.1, exactly as in the proof of Theorem 4.2
we arrive at
$$\left |U_{n}(f g)(z)-U_{n}(f)(z)U_{n}(g)(z)-\frac{2z(1-z)f^{\prime}(z)g^{\prime}(z)}{n}\right |\le \frac{C(r, f, g)}{n^{2}},$$
for all $n\in \mathbb{N}$ and $|z|\le r$, with $C(r, f, g)>0$ independent of $n$ and depending on $r, f, g$. $\hfill \square$

Based on the estimate in Theorem 5.2 and following similar reasonings with those in the proofs of Corollaries 4.3 and 4.4, we easily arrive at the next result.

{\bf Corollary 5.3.} {\it Suppose that $R > r \ge 1$ and $f, g:\mathbb{D}_{R}\to \mathbb{C}$ are analytic in $\mathbb{D}_{R}=\{z\in \mathbb{C} ; |z|<R\}$, that is $f(z)=\sum_{k=0}^{\infty}a_{k}z^{k}$ and $f(z)=\sum_{k=0}^{\infty}b_{k}z^{k}$ for all $z\in \mathbb{D}_{R}$. If $f$ and $g$ are not constant functions, then for any $1\le r <R$ we have
$$\|U_{n}(f g)-U_{n}(f)U_{n}(g)\|_{r}\sim \frac{1}{n}, \, n\in \mathbb{N},$$
where the constants in the equivalence are independent of $n$ but depend on $r, f, g$.}

\section{Results for Bernstein-Faber operators}

In this section we extend the results in Section  4 to Bernstein-Faber polynomials attached to compact sets of the complex plane.
For this purpose, firstly let us briefly recall some classical concepts and results about Faber polynomials (for more details see \cite{Gaier}, \cite{Suetin}).

If $G\subset \mathbb{C}$ is a compact set such that $\tilde{\mathbb{C}}\setminus G$ is connected, denote by $A(G)$ the Banach space of all functions that are continuous on $G$ and analytic in the interior of $G$, endowed with the uniform norm $\|f\|_{G}=\sup\{|f(z)| ; z\in G\}$. If we denote $\mathbb{D}_{r}=\{z\in \mathbb{C} ; |z|<r\}$ then according to the Riemann Mapping Theorem, a unique
conformal mapping $\Psi$ of $\tilde{\mathbb{C}}\setminus \overline{\mathbb{D}}_{1}$ onto $\tilde{\mathbb{C}}\setminus G$ exists so that $\Psi(\infty)=\infty$ and $\Psi^{\prime}(\infty)>0$. The $n$-th Faber polynomial $F_{n}(z)$ attached to $G$ may be defined by
$$\frac{\Psi^{\prime}(w)}{\Psi(w)-z}=\sum_{n=0}^{\infty}\frac{F_{n}(z)}{w^{n+1}}, \, \, z\in G, |w|>1.$$
Then $F_{n}(z)$ is a polynomial of exact degree $n$.

If $f\in A(G)$ then
$$a_{n}(f)=\frac{1}{2\pi i}\int_{|u|=1}\frac{f[\Psi(u)]}{u^{n+1}}du=\frac{1}{2\pi }\int_{-\pi}^{\pi}f[\Psi(e^{it})]e^{-int}dt, n\in \mathbb{N}\cup \{0\}$$
are called the Faber coefficients of $f$ and $\sum_{n=0}^{\infty}a_{n}(f)F_{n}(z)$ is called the Faber expansion (series)
attached to $f$ on $G$. (Here $i^{2}=-1$.) The Faber series represent a natural generalization of Taylor series
when the unit disk is replaced by an arbitrary simply connected domain bounded by a "nice" curve.

In \cite{Gal}, p. 19, the Bernstein-Faber polynomials were defined by the formula
$${\cal{B}}_{n}(f ; G)(z)=\sum_{p=0}^{n}{n \choose p}\Delta_{1/n}^{p}F(0)\cdot F_{p}(z), z\in G, \ n\in \mathbb{N},$$
where
$$\Delta^{p}_{h}F(0)=\sum_{k=0}^{p}(-1)^{p-k}{p \choose k}F(kh), \, F(w)=\frac{1}{2\pi i}\int_{|u|=1}\frac{f(\Psi(u))}{u-w}du, \ w\in \mathbb{D}_{1}.$$ Here, since $F(1)$ is involved in $\Delta^{n}_{1/n}F(0)$ and in the definition of ${\mathcal{B}}_{n}(f ; G)(z)$ too, in addition we will suppose that $F$ can be extended by continuity on the boundary $\partial \mathbb{D}_{1}$.

{\bf Remarks.} 1) For $G=\overline{\mathbb{D}}_{1}$, since $\Psi(z)=z$, $F(z)=f(z)$ and $F_{p}(z)=z^{p}$, it is easy to see that the above Bernstein-Faber polynomials one reduce to the classical complex Bernstein polynomials given by
$$B_{n}(f)(z)=\sum_{p=0}^{n}{n \choose p}\Delta_{1/n}^{p}f(0)z^{p}=\sum_{p=0}^{n}{n \choose p}z^{p}(1-z)^{n-p}f(p/n).$$

2) It is known that, for example, $\int_{0}^{1}\frac{\omega_{p}(f\circ \Psi ; u)_{\partial \mathbb{D}_{1}}}{u}du<\infty$ is
a sufficient condition for the continuity on $\partial \mathbb{D}_{1}$  of $F$ in the above
definition of the Bernstein-Faber polynomials (see e.g. \cite{Gaier}, p. 52, Theorem 6). Here $p\in \mathbb{N}$ is arbitrary fixed.

Now, we are in position to prove the extensions of the results in Section 4, as follows.

{\bf Theorem 6.1.} {\it Let $G$ be a continuum (that is a
connected compact subset of $\mathbb{C}$) and suppose that $f, g$ are
analytic in $G$, that is there exists $R>1$ such that $f$ and $g$ are
analytic in $G_{R}$, that is $f(z)=\sum_{k=0}^{\infty}a_{k}(f)F_{k}(z)$ and $g(z)=\sum_{k=0}^{\infty}a_{k}(g)F_{k}(z)$, for all $z\in G_{R}$. Here recall that $G_{R}$ denotes the interior
of the closed level curve $\Gamma_{R}$ given by $\Gamma_{R}=\{z ;
|\Phi(z)|=R \} =\{\Psi(w) ; |w|=R \}$ (and that $G\subset
\overline{G}_{r}$ for all $1<r<R$). Also, we suppose that $F$
given in the definition of Bernstein-Faber polynomials can be
extended by continuity on $\partial \mathbb{D}_{1}$.

Let $1<r <R$.

(i) (Gr\"uss estimate) For all $z\in \overline{G_{r}}$ and $n\in \mathbb{N}$ we have
$$|{\cal{B}}_{n}(f g ; G)(z)-{\cal{B}}_{n}(f ; G)(z)\cdot {\cal{B}}_{n}(g; G)(z)|\le \frac{C}{n},$$
where $C$ depends on $f$, $g$ and $r$ but is independent of $n$.

(ii) (Gr\"uss-Voronovskaya estimate) For all $z\in \overline{G_{r}}$ and $n\in \mathbb{N}$ we have
$$\left |{\cal{B}}_{n}(f g ; G)(z)-{\cal{B}}_{n}(f ; G)(z)\cdot {\cal{B}}_{n}(g ; G)(z)\right .$$
$$\left .-\sum_{k=2}^{\infty}\frac{k(k-1)}{2 n}[F_{k-1}(z)-F_{k}(z)]\cdot [a_{k}(f g)-f(z)a_{k}(g)-g(z)a_{k}(f)]\right |\le \frac{C}{n^{2}},$$
where $C$ depends on $f, g, r$ but is independent of $n$.}

{\bf Proof.} (i) We can write
$$|{\cal{B}}_{n}(f g ; G)(z)-{\cal{B}}_{n}(f ; G)(z)\cdot {\cal{B}}_{n}(g; G)(z)|
\le |{\cal{B}}_{n}(f g ; G)(z)-f(z)g(z)|$$
$$+|f(z)|\cdot |g(z)-{\cal{B}}_{n}(g ; G)(z)|+|{\cal{B}}_{n}(g ; G)(z)|\cdot |f(z)-{\cal{B}}_{n}(f ; G)(z)|$$
$$\le \frac{C(f, g)}{n}+\|f\|_{\overline{G_{r}}}\cdot \frac{C(g)}{n}+ M(g)\cdot \frac{C(f)}{n}=\frac{C}{n},$$
where we used the result in \cite{Gal}, p. 20, which states that
$$|{\cal{B}}_{n}(f; G)(z)-f(z)|\le \frac{C(f)}{n}, \mbox{ for all } z\in \overline{G_{r}} \mbox{ and } n\in \mathbb{N}.$$

(ii) Firstly, we have
$${\cal{B}}_{n}(f g ; G)(z)-{\cal{B}}_{n}(f ; G)(z)\cdot {\cal{B}}_{n}(g ; G)(z)$$
$$-\sum_{k=2}^{\infty}\frac{k(k-1)}{2 n}[F_{k-1}(z)-F_{k}(z)]\cdot [a_{k}(f g)-f(z)a_{k}(g)-g(z)a_{k}(f)]$$
$$=\left [{\cal{B}}_{n}(f g ; G)(z)-(f g)(z)-\sum_{k=2}\frac{k(k-1)}{2 n}\cdot a_{k}(f g)[F_{k-1}(z)-F_{k}(z)]\right ]$$
$$-f(z)\left [{\cal{B}}_{n}(g ; G)(z)-g(z)- \sum_{k=2}\frac{k(k-1)}{2 n}\cdot a_{k}(g)[F_{k-1}(z)-F_{k}(z)]\right ]$$
$$-g(z)\left [{\cal{B}}_{n}(f ; G)(z)-f(z)-\sum_{k=2}\frac{k(k-1)}{2 n}\cdot a_{k}(f)[F_{k-1}(z)-F_{k}(z)]\right ]$$
$$+[g(z)-{\cal{B}}_{n}(g : G)(z)]\cdot [{\cal{B}}_{n}(f : G)(z)-f(z)].$$
Then, taking into account the estimate in \cite{Gal}, p. 20 mentioned at the above point (i) and the Voronovskaya-type estimate in \cite{Gal2}, p. 88,
Theorem 1.11.4, (i), given by
$$|{\cal{B}}_{n}(f ; G)(z)-f(z)-\sum_{k=2}\frac{k(k-1)}{2 n}\cdot a_{k}(f)[F_{k-1}(z)-F_{k}(z)|\le \frac{C}{n^{2}},$$
the proof is immediate. $\hfill \square$.

{\bf Remarks.} 1) When $G=\mathbb{D}_{R}$ then Theorem 6.1, (ii) one reduces to Theorem 4.2.

2) The are many concrete examples for $G$ when the conformal mapping $\Psi$ and the Faber polynomials associated to $G$ (and consequently when the Bernstein-Faber polynomials too) can explicitly be calculated (see for details \cite{Gal2}, pp. 81-83) :

(i) $G$ is the continuum bounded by the $m$-cusped hypocycloid $H_{m}$ ($m=2, 3, ...,$), given by the parametric equation
$$z=e^{i \theta}+\frac{1}{m-1}e^{-(m-1)i\theta}, \theta\in [0, 2\pi),$$
case when $\Psi(w)=w+\frac{1}{(m-1)w^{m-1}}$ and the Faber polynomials can explicitly be calculated ;

(ii) a) $G$ is the regular $m$-star ($m=2, 3, ...,$) given by
$$S_{m}=\{x \omega^{k} ; 0\le x \le 4^{1/m}, k=0, 1, ..., m-1, \, \omega^{m}=1\},$$
case when $\Psi(w)=w\left (1+\frac{1}{w^{m}}\right )^{2/m}$ and the Faber polynomials can explicitly be calculated ;

(iii) $G$ is the $m$-leafed symmetric lemniscate, $m=2, 3, ..., $ with its boundary given by
$$L_{m}=\{z \in \mathbb{C} ; |z^{m}-1|=1\},$$
case when $\Psi(w)=w \left (1+\frac{1}{w^{m}}\right )^{1/m}$  and the Faber polynomials can explicitly be calculated ;

(iv) $G$ is the semidisk
$$SD=\{z\in \mathbb{C} ; |z|\le 1 \mbox{ and } |Arg(z)|\le \pi/2\},$$
case when $\Psi(w)=\frac{2(w^{3}-1)+3(w^{2}-w)+2(w^{2}+w+1)^{3/2}}{w(w+1)\sqrt{3}}$
and the attached Faber polynomials can explicitly  be calculated  ;

(v) $G$ is a circular lune or $G$ is an annulus sector, cases when again the conformal mapping $\Psi$ and the Faber polynomials can explicitly be calculated.

\end{document}